\documentclass[12pt,a4paper]{article}

\usepackage{amsmath}
\usepackage{amssymb}

 \newcommand{\NN}{\mathbb N}
 \newcommand{\RR}{\mathbb R}
 \newcommand{\CC}{\mathbb C}
 \newcommand{\ZZ}{\mathbb Z}

 \newcommand{\gG}{\varGamma}
 \newcommand{\gP}{\varPhi}

 \newtheorem{theorem}{Theorem}
 \newtheorem{lemma}{Lemma}
 \newtheorem{prop}{Proposition}

 \newcommand{\qed}{\hfill $\square$}
 \newcommand{\dens}{{\rm dens}}

\begin{document}

\centerline{\Large \bf Diffraction of weighted lattice subsets}

\vspace{15mm}
        
\centerline{{\sc Michael Baake}}
\vspace{5mm} 

{\small \noindent
\begin{center}
Institut f\"ur Mathematik, Universit\"at Greifswald \\
Jahnstr.~15\thinspace a, 17487 Greifswald, Germany 
\end{center} }

\vspace{10mm}

\centerline{Dedicated to Robert V.\ Moody on the occasion 
          of his 60th birthday}

\vspace{10mm}
\begin{abstract} 
A Dirac comb of point measures in Euclidean space with bounded complex 
weights that is supported  on a lattice $\gG$\/ inherits certain general 
properties from the lattice structure.
In particular, its autocorrelation admits a factorization into a continuous
function and the uniform lattice Dirac comb, and its diffraction measure
is periodic, with the dual lattice $\gG^*$ as lattice of periods. This
statement remains true in the setting of a locally compact Abelian
group that is also $\sigma$-compact.

\end{abstract}

\bigskip
\bigskip

\section{\Large Introduction}

Mathematical diffraction theory is concerned with the Fourier analysis
of the autocorrelation measures of unbounded, but translation bounded,
complex measures in Euclidean space, compare \cite{Cowley,Wel,Hof,BH,Lagarias} 
and references listed there. More generally, the same question is also analyzed
in the setting of locally compact Abelian groups \cite{AL,BF,Sch}.
There are many interesting and relevant open questions connected with it, 
in particular how to assess the spectral type of the (positive) diffraction 
measure, which is the Fourier transform of the autocorrelation.

One important subclass of translation bounded measures
consists of the set of discrete Dirac combs which
are supported on a lattice or a subset thereof. 
This paper provides a characterization of the autocorrelation and
diffraction measures associated to measures of this class. As an
application, we derive a relation between the diffraction measures
of complementary lattice subsets. Other recent (sometimes implicit)
applications of this characterization include deterministic cases such as
the visible lattice points or the $k$th power free integers \cite{BMP,BMnew}
and the class of substitution lattice systems \cite{Lee,LMS},
but also a wide class of lattice gas models \cite{BH,Hoeffe}. The latter are
presently gaining practical importance in crystallography due to the need for 
a better understanding of diffuse scattering, compare \cite{BH,BH2,Wel} and 
references given there. 

The article is organized as follows. First, we describe the problem in the
Euclidean setting and state the corresponding result in Theorem~\ref{thm1}.
The next Section then covers the step by step proof of this result.
All details, or at least precise references, are given in order to make
the presentation relatively self-contained and accessible also to readers
with a more applied background. As general references for topological
concepts and results, we use \cite{Bou,Q}.
As one application, we then compare lattice subsets with their
complements, with special emphasis on the homometry problem
(see Theorem~\ref{complement}). This is followed by an extension of 
Theorem~\ref{thm1} to the more general
situation that Euclidean space $\RR^n$ is replaced by a locally compact
Abelian group (LCA group) $G$ which is Polish. An analogous result is true here, 
and is summarized in Theorem~\ref{thm2}. Finally, we close with a slightly 
weaker result for the situation of $\sigma$-compact LCA groups
in Theorem~\ref{thm4}.

\section{\Large Euclidean lattices}

Let $\gG$ be a lattice in $\RR^n$, i.e.\ the integer span of $n$
vectors that are linearly independent over $\RR$. We are interested
in the {\em weighted\/} Dirac comb \cite{Cordoba}
\begin{equation} \label{comb1}
   \omega \; = \; \sum_{t\in\gG} w(t)\, \delta_t
\end{equation}
where $w\!:\gG\rightarrow\CC$ is a bounded function and
$\delta_t$ denotes the normalized point (or Dirac) measure located at 
$t$, so that $\delta_t(\varphi) = \varphi(t)$ for continuous functions
$\varphi$. Since $\gG$ is uniformly discrete and $w$ is bounded, $\omega$ 
is a translation bounded (or shift bounded) complex measure, cf.\
\cite[Ch.\ I.1]{BF}. By a measure, we always mean a regular (complex) 
Borel measure, identified with the corresponding linear functional on
${\mathcal K}(\RR^n)$, the space of continuous functions of compact
support. This is justified by the Riesz-Markov theorem, compare
\cite{RS,BF,Hof,BH} for details. The space of measures is then equipped
with the vague topology which we will use throughout.

Let us first assume that the {\em natural autocorrelation measure\/} of 
$\omega$ exists. It is defined with respect to an increasing sequence of balls
around the origin, and denoted by $\gamma^{}_{\omega}$. So, we assume for
the moment that
\begin{equation} \label{auto1}
    \gamma^{}_{\omega} \; = \; \lim_{r\to\infty}
    \frac{\omega^{}_r * \tilde{\omega}^{}_r}{{\rm vol}(B_r(0))}  
\end{equation}
exists as a vague limit, where $\omega^{}_r$ means the restriction of 
$\omega$ to the (open) ball $B_r(0)$ of radius $r$ around $0$,
$\tilde{\omega}^{}_r := (\omega^{}_r)\tilde{\hphantom{t}}$,
and $\tilde{\omega}$ is the measure defined by 
$\tilde{\omega} (\varphi) = \overline{\omega (\tilde{\varphi})}$
for any continuous function $\varphi$ of compact support,
with $\tilde{\varphi}(x) := \overline{\varphi(-x)}$. In particular,
for $\omega$ of (\ref{comb1}), we have $\tilde{\omega} =
\sum_{t\in\gG} \overline{w(t)}\, \delta_{-t}$.

Since $\gG$ is a lattice, the natural autocorrelation is a pure point
measure of the form
$\gamma^{}_{\omega} = \sum_{z\in\gG} \nu(z)\, \delta_z$
with the coefficients being given by
\begin{eqnarray} \label{auto2}
    \nu(z) & = & \lim_{r\to\infty} \,
    \frac{1}{{\rm vol}(B_r(0))} \, 
    \sum_{\substack{t,t'\in\gG_r \\ t - t' = z}}
    {w(t)}\, \overline{w(t')}  \\
    & = &   \lim_{r\to\infty} \,
    \frac{1}{{\rm vol}(B_r(0))} \, \sum_{t\in\gG_r}
    {w(t)}\, \overline{w(t-z)}  \nonumber
\end{eqnarray}
Here, $\gG_r = \gG\cap B_r(0)$, and there are several
variants to write $\nu(z)$ as a limit. For fixed $z$ and $r \gg |z|$, the 
two approximations for finite radius used in (\ref{auto2}), 
prior to dividing by ${\rm vol}(B_r(0))$, 
differ by surface contributions which are uniformly small in comparison to 
the bulk. To make this precise, one needs an estimate for the
number of lattice points in spherical shells of thickness $s$ and
radius $r \gg s$ which is $O(r^{n-1})$ in $n$ dimensions, see \cite{KR}.
Hence, the above limit does not depend on such details, and we
will make use of this freedom in the Euclidean setting without further 
mentioning.

Since $\gG$ is a lattice in $\RR^n$, its density is well-defined as a 
limit, namely
\begin{equation} \label{density}
   {\rm dens}(\gG) \; = \; \lim_{r\to\infty} \,
   \frac{|\gG\cap B_r(a)|}{{\rm vol}(B_r(0))}
\end{equation}
with $a\in\RR^n$ and $|A|$ denoting the cardinality of a (finite) set $A$. 
This limit exists {\em uniformly\/} in $a$, see \cite{Hof,Sch}. 
As a lattice, $\gG$ is uniformly discrete (i.e.\ the minimal distance
between any two distinct points is strictly positive), and the natural
autocorrelation thus exists if and only if the coefficients $\nu(z)$ of
(\ref{auto2}) exist for all $z\in\gG$.

So far, we explicitly assumed
the existence of the natural autocorrelation as a limit, and hence
its uniqueness. This is the standard situation investigated e.g.\ in
crystallography, and it seems very adequate in view of the usual
homogeneity of the systems analyzed. Mathematically, however, this
assumption is not necessary, and we will thus formulate our results
in a slightly more general setting.

Our key assumption is that the complex weight function 
$w$ is bounded, so that $\omega$ is translation bounded. This
is already sufficient to ensure that all finite approximations to the
autocorrelation are uniformly translation bounded, 
see \cite[Prop.\ 2.2]{Hof}. Consequently, there is always at
least one vague limit point, but, in general, we might have
several. To make this precise, let us define
$$ \gamma^{(r)}_{\omega} \; = \; 
   \frac{\omega^{}_r * \tilde{\omega}^{}_r}{{\rm vol}(B_r(0))}\; $$
and consider the family of autocorrelation approximants
$\{\gamma^{(r)}_{\omega} \mid r > 0 \}$, which are all positive
definite measures by construction. When restricted to a compact
set $K$, the corresponding family of finite measures is precompact
in the vague topology.

Let us now consider {\em any\/} vague limit point $\gamma$ of the family
$\{\gamma^{(r)}_{\omega} \mid r > 0 \}$. It is then possible to select a 
sequence of radii $(r^{}_i)^{}_{i\in\NN}$ such that 
$\gamma^{(r^{}_i)}_{\omega}\to\gamma$ vaguely as $i\to\infty$.
As we will see, the proof of our main result (Theorem~\ref{thm1})
will only depend on the existence of such a sequence, so we can
formulate it as a result for each limit point.

The standard tools from Fourier analysis are also needed. The Fourier 
transform of a rapidly decreasing (or Schwartz) function $\varphi$ is
$$ \hat{\varphi}(x) \; = \; \int_{\RR^n} e^{- 2\pi i x y} \, \varphi(y)
   \, {\rm d}y  $$
where $x y$ is the standard Euclidean scalar product in $\RR^n$.
From here, we take the usual route to the Fourier transform of measures
that are, at the same time, tempered distributions, see \cite[Ch.\ IX.1]{RS} 
or \cite[Ch.\ 7]{Rudin2} for general background and \cite{BH} for our 
conventions. An alternative approach, without any reference to Schwartz
functions, could follow \cite[Ch.\ I.4]{BF} and employ the positive
definiteness of $\gamma^{}_{\omega}$ which guarantees the existence
of $\hat{\gamma}^{}_{\omega}$.

We call a Borel set $A$ {\em regular\/}, if $\partial A$ has Lebesgue
measure $0$, and say that a measure $\mu$ is {\em supported\/} on a regular
Borel set $A$, if $\mu$ and $\mu|_A$, the restriction of $\mu$ to $A$,
are identical as measures on $\RR^n$.
The central result can now be phrased as follows.
\begin{theorem} \label{thm1}
   Let $\gG$ be a lattice in $\RR^n$ and $\omega$ a weighted Dirac comb 
   on $\gG$ with bounded complex weights. Let $\gamma^{}_{\omega}$ 
   be any of its autocorrelations, i.e.\ any of the limit points of
   the family  $\{\gamma^{(r)}_{\omega} \mid r > 0 \}$.
   Then the following holds.

$(1)$ The autocorrelation $\gamma^{}_{\omega}$ can be represented as
$$ \gamma^{}_{\omega} \; = \; \gP \cdot \delta^{}_{\gG} 
   \; = \; \sum_{x\in\gG} \gP(x)\,\delta^{}_x $$
   where $\gP\!:\RR^n\to\CC$ is a bounded continuous 
   positive definite function that 
   interpolates the autocorrelation coefficients $\nu(x)$ as defined 
   at $x\in\gG$. Moreover, there exists 
   such a\/ $\gP$ which extends to an entire function $\gP\!:\CC^n\to\CC$
   with the additional growth restriction that there are constants
   $C,R \ge 0$ and $N\in\ZZ$ such that
   $|\gP(z)| \le C (1+|z|)^N \exp(R\, |{\rm Im}(z)|)$
   for all $z\in\CC^n$.  

$(2)$ The Fourier transform $\hat{\gamma}^{}_{\omega}$, also called a
   diffraction measure of $\omega$, is a translation bounded positive 
   measure that is periodic with the dual lattice 
   $\gG^* = \{u\mid uv\in\ZZ \,\mbox{ for all }\, v\in\gG\}$ as lattice 
   of periods. Furthermore,
   $\hat{\gamma}^{}_{\omega}$ has a representation as a convolution,
$$ \hat{\gamma}^{}_{\omega} \; = \; \varrho * \delta^{}_{\gG^*}\, , $$
   in which $\varrho$ is a finite positive measure supported on a 
   fundamental domain of $\gG^*$ that is contained in the ball
   of radius $R$ around the origin.
\end{theorem}

\noindent
{\sc Remarks}: 
The interpolating function $\gP$ is not unique. It can be changed by
adding any positive definite continuous (entire) function which
vanishes on $\gG$. Also, the positive measure $\varrho$ depends
on the choice of the fundamental domain of $\gG^*$. As soon as the
latter is fixed (as a regular Borel set, say, to avoid pathologies
with singular versus continuous parts), $\varrho$ is unique. This is
so because $\hat{\gamma}^{}_{\omega}$ is $\gG^*$-periodic and thus
defines a unique measure on the factor group $\RR^n/\gG^*$, which,
in turn, gives rise to a unique measure $\varrho$ on the fundamental
domain chosen.

A rather natural fundamental domain of $\gG^*$ can be constructed
from the Voronoi cell of $\gG^*$, which is a polytope, namely the set
of all points of $\RR^n$ which are not farther apart from the origin
than from any other point of $\gG^*$. It is the intersection of
finitely many closed half-spaces and hence a zonotope.
A fundamental domain emerges from it by removal of certain boundaries,
starting with the $(n-1)$-boundaries, then the $(n-2)$-boundaries,
and so on down to the 0-boundaries. Due to the lattice structure of
$\gG^*$, these boundaries always come in translation pairs,
one member of which is removed.
What remains, is a regular Borel set $E$ which is a true fundamental domain of 
$\gG^*$, i.e.\ the $\gG^*$-translates of $E$ form a partition of $\RR^n$.
Its existence is vital
for the above statement because $\varrho$ may contain point measures.
 
The radius $R$ in the above growth estimate can be taken as the radius 
of the circumsphere of the fundamental domain of $\gG^*$ chosen. 
This follows from the Paley-Wiener 
theorem for distributions with compact support \cite[Thm.\ 7.23]{Rudin2}.
If a fundamental domain is constructed on the basis of the Voronoi cell of 
$\gG^*$, then its circumradius $R$, the covering radius of $\gG^*$, is minimal
among all fundmanental domains of $\gG^*$.

If there is only one limit point to the family 
$\{\gamma^{(r)}_{\omega} \mid r > 0 \}$, we are in the standard situation
of crystallography, where the homogeneity of the systems makes this
a very natural assumption. In this case, Theorem~\ref{thm1} simply
refers to {\em the\/} autocorrelation.

The above mentioned question of the spectral type 
of $\hat{\gamma}^{}_{\omega}$ is now obviously reduced to that of 
the spectral type of $\varrho$. The latter is 
a finite positive measure and admits the unique decomposition
$$ \varrho \; = \; (\varrho)^{}_{pp} 
        + (\varrho)^{}_{sc} + (\varrho)^{}_{ac} $$
with respect to Lebesgue measure, which is the natural reference
measure in this context. As usual, $pp$, $sc$ and $ac$ stand for
pure point, singular continuous and absolutely continuous, see
\cite[Sec.\ I.4]{RS} for details. 

The result of Theorem~\ref{thm1} is not really restricted to Dirac combs.
If $h$ is a (continuous) $L^1$-function, say, and $\omega$
a Dirac comb of the above type, then $h*\omega$ is a well-defined
translation bounded measure, with autocorrelation
$\gamma = (h*\tilde{h}) * \gamma^{}_{\omega}$ and diffraction measure
$\hat{\gamma} = |\hat{h}|^2 \, \hat{\gamma}^{}_{\omega}$ by the
convolution theorem. This can be considered as a situation where
$h$ describes a more realistic profile of the scatterers (e.g.\ atoms),
and Theorem \ref{thm1} can then be applied to $\gamma^{}_{\omega}$
and $\hat{\gamma}^{}_{\omega}$. Convolution is also a key ingredient to
tackle the problem of diffraction at high temperature, compare
\cite{Hof-HT}.

Various applications were already mentioned in the Introduction.
Also, the case that $\varrho  =  (\varrho)^{}_{pp}$ can now be analyzed
in more detail, and a set of necessary and sufficient conditions will
be given in \cite{BMnew}. Arrangements of scatterers where
$\hat{\gamma}^{}_{\omega}$ is a pure point measure are often
called pure point diffractive. These cases are of particular interest
because they include perfect crystals, but also several aperiodic lattice
substitution systems, see \cite{BMS,Lee} for details.

\section{\Large Proof of Theorem 1} \label{proofs}

The key idea is to use an appropriate regularization of the point measure 
$\omega$ of (\ref{comb1}) for the construction of the function $\gP$ and the
measure $\varrho$. It relies on the existence of a properly converging
sequence of approximating measures.
We will employ Ascoli's Theorem to construct a (Lipschitz)
continuous interpolation of the autocorrelation coefficients of (\ref{auto2}), 
and a combination of Bochner's Theorem with the convolution property of the
Fourier transform to conclude on the existence of the measure $\varrho$.

Let $c\in C^{\infty}_c$ be a non-negative
bump function with compact support contained in $B_{\varepsilon}(0)$
where $\varepsilon > 0$ is at most half the packing radius of $\gG$
(the technical reason for this restriction will become clear below). 
Let $c(0)=c^{}_0$ be the maximal value of $c$, so that 
$0\le c(x)\le c^{}_0$ for all $x\in\RR^n$. 
Such a function, which is also called an approximate identity, can be 
constructed as follows, compare \cite[p.\ 168]{Lang}. Start from
$\phi(x) = \exp\big(|x|^2/(|x|^2 - 1)\big)$ for $|x|<1$ and $\phi(x)=0$
otherwise. This is a $C^{\infty}$-function with support in the closed
unit ball, with $\phi(0)=1$. Now, set $c(x) = c^{}_0\, \phi(x/\varepsilon)$.
Clearly, $0\le c(x)\le c^{}_0$, the support is in the closed ball of radius
$\varepsilon$, and $c$ is integrable, i.e.\ $\|c\|^{}_1 < \infty$.
Also, $(c*\tilde{c}) (0) = \int_{\RR^n} c(x) \tilde{c}(-x)\,{\rm d}x
= \|c\|^2_2 < \infty$, and, for fixed $\varepsilon > 0$, 
we can always adjust $c^{}_0$ so that
$\|c\|^2_2 = 1$. We henceforth assume that such a function $c$ is given.

Recall that a function $c$ is (globally) Lipschitz, with Lipschitz constant
$L_c < \infty$, if $|c(x) - c(y)| \le L_c\, |x-y|$ for all $x,y\in\RR^n$, and
$L_c$ is the smallest number with this property.
In particular, $c$ is then uniformly continuous.
\begin{lemma} \label{lip}
   A non-negative bump function $c\in C^{\infty}_c$ is 
   Lipschitz. Furthermore, also $\tilde{c}$ and $c * \tilde{c}$ are
   Lipschitz, with $L^{}_{\tilde{c}} = L^{}_c$ and
   $L^{}_{c*\tilde{c}} \le \|c\|^{}_1 \, L^{}_c$.
\end{lemma}
{\sc Proof}: It is clear that $c$ is Lipschitz because it is smooth and 
has compact support, so 
$L_c \le \sup_{x\in\RR^n} \big|{\rm grad} ( c(x)) \big| < \infty$. 
Since $c$ is a real function, we
have $\tilde{c}(x)=\overline{c(-x)}=c(-x)$, so that $L^{}_c = L^{}_{\tilde{c}}$
is obvious. Finally, consider
\begin{eqnarray*}
  \big| (c*\tilde{c})(x) - (c*\tilde{c})(y) \big| & = &
  \bigg| \int_{\RR^n} c(z)\, \big(\tilde{c}(x-z) - \tilde{c}(y-z)\big)
        \, {\rm d}z \,\bigg|  \\
  & \le & \int_{\RR^n} c(z)\, \big|\tilde{c}(x-z) - \tilde{c}(y-z)\big|
        \, {\rm d}z  \\
  & \le & L_c\, |x-y| \int_{\RR^n} c(z) \, {\rm d}z
  \;\, = \;\, \|c\|^{}_1 \, L_c\, |x-y|
\end{eqnarray*}
which establishes the assertion. \qed

\smallskip
Fix a radius $r>0$ and set $\gG_r = \gG\cap B_r(0)$. Let
$\omega^{}_r = \sum_{t\in\gG_r} w(t) \delta_t$ be the weighted Dirac 
comb of $\gG_r$. Define $f_r = c * \omega^{}_r$, so that
\begin{equation} \label{cont-fun}
   f_r(z) \; = \; \sum_{t\in\gG_r} w(t)\, c(z-t) \, .
\end{equation} 
This is a continuous function,
as is $\tilde{f}_r$. A simple calculation then shows that
$ f_r * \tilde{f}_r = (c * \tilde{c}) * (\omega^{}_r * \tilde{\omega}^{}_r)$
which is the continuous positive definite function
\begin{equation} \label{faltung}
   \big(f_r * \tilde{f}_r\big) (z) \; = \; \sum_{u,v\in\gG_r}
    w(u)\, \overline{w(v)}\, \big(c * \tilde{c}\big) (z - u - v)\, . 
\end{equation}
Due to the small support of $c$ (recall that $\varepsilon$ is at most
half the packing radius of $\gG$), and hence that of $c * \tilde{c}$,
most terms of the double sum vanish, and it effectively
collapses to a single sum. In particular, if $z=t$ is in the lattice,
$\big(c * \tilde{c}\big) (t-u-v)$ is either 0 or 1 by construction. In fact,
for $u\in\gG$, the value 1 is taken precisely for $v=t-u$.
This is a consequence of the above choice of $\varepsilon$, and
also motivates it, a posteriori.

Let us now define $g^{}_r = \frac{1}{{\rm vol}(B_r(0))}\, f_r * \tilde{f}_r$.
Then it is clear from the above that $g^{}_r(t)\to\nu(t)$ as $r\to\infty$
for all $t\in\gG$. We have thus constructed a family of continuous
positive definite functions $\{g^{}_r\mid r>0\}$ that pointwise converges to 
the autocorrelation coefficients for all $t\in\gG$.
We now have to check whether it will result in a continuous interpolation.
\begin{lemma} \label{g-lip}
   The family $\{g^{}_r\mid r>0\}$ is uniformly Lipschitz, hence in
   particular equi-continuous, and uniformly bounded.
\end{lemma}
{\sc Proof}: 
The function $f_r$ of (\ref{cont-fun}), for any $r>0$, is a 
finite linear combination of Lipschitz functions and hence Lipschitz,
by Lemma \ref{lip}.
With $W:=\sup_{t\in\gG} |w(t)|$, we can estimate the Lipschitz
constant of $f_r$ as $L^{}_{f^{}_r}\le |\gG_r|\, W\, L^{}_c$ since
the support of $c$ is so small that no two bumps in $f_r$ overlap. 
Lemma~\ref{lip} then gives $L^{}_{f^{}_r * \tilde{f}^{}_r} \le
\|f_r\|^{}_1 \, L^{}_{f^{}_r}$, and with $\|f_r\|^{}_1 \le W\, \|c\|^{}_1$
we obtain
$$ L^{}_{f^{}_r * \tilde{f}^{}_r} \; \le \; 
   |\gG_r| \, W^2 \, \|c\|^{}_1 \, L^{}_c  \; < \; \infty \, .$$
Consequently, also $g^{}_r$ is Lipschitz, and we have
$$ L^{}_{g^{}_r} \; \le \; \frac{|\gG_r|}{{\rm vol}(B_r(0))}\;
    W^2 \, \|c\|^{}_1 \, L_c \, .$$
The first factor on the right hand side is the number of
lattice points in a ball of radius $r$ divided by its volume.
This is known \cite{KR} to be
\begin{equation} \label{spheres}
   \frac{|\gG_r|}{{\rm vol}(B_r(0))} \; = \;
    {\rm dens}(\gG) + {\mathcal O} (r^{-1})
    \qquad \mbox{as $r\to\infty$} 
\end{equation}
from which the uniform Lipschitz condition follows, and hence also
the equi-continuity of the family $\{ g^{}_r \mid r>0 \}$.

Similarly, consider $g^{}_r = (c * \tilde{c}) * 
\frac{\omega^{}_r * \tilde{\omega}^{}_r}{{\rm vol}(B_r(0))}$ which is a
continuous function composed of non-overlapping little bumps.
Consequently, using (\ref{auto1}) and (\ref{auto2}) and observing that
$\|c * \tilde{c}\|^{}_{\infty} = \|c\|^2_2 = 1$ (due to our choice of
$c^{}_0$), we see that
$$ \big| g^{}_r (z) \big| \; \le \; \big\| g^{}_r \big\|_{\infty} \; \le \;
    \frac{|\gG_r|}{{\rm vol}(B_r(0))}\, W^2  $$
from which, again with (\ref{spheres}), one obtains equi-boundedness.
\qed

\smallskip
Lemma \ref{g-lip} allows us to use Ascoli's Theorem, see 
\cite[Cor.\ III.3.3]{Lang} for a version that covers our situation
(note that $\RR^n$ is only locally compact, but $\sigma$-compact, i.e.\
it can be covered by a countable family of compact sets).
In particular, no matter whether we start from the entire family 
$\{g^{}_r \mid r>0\}$ or from a sequence
$\big( g^{}_{r_i} \big)_{i\in\NN}$ (as needed for the case that we
have several autocorrelations), there is
a (sub)sequence of radii, $(r_j)^{}_{j\in\NN}$, such that, as $j\to\infty$,
the $g^{}_{r^{}_j}$ converge compactly to some function $g$ which is then 
bounded and continuous. Since the family is actually equi-uniformly
continuous, convergence is globally uniform.
As a limit of positive definite functions, $g$ is still positive definite. 
Moreover, it is also (globally) Lipschitz, and $g(t) = \nu(t)$ for all 
$t\in\gG$ by construction. So far, we have established:
\begin{prop} \label{rep1}
   Let $\gG$ be a lattice in $\RR^n$ and\/ $\omega$ the Dirac comb of 
   $(\ref{comb1})$, with bounded complex weights.
   Let $\gamma^{}_{\omega}$ be any of its $($natural\/$)$ 
   autocorrelation measures, the latter defined as the limit 
   of a suitable sequence $\big( \gamma^{(r_i)}_{\omega}\big)_{i\in\NN}$ of
   approximants, with $(r_i)^{}_{i\in\NN}$ an increasing, unbounded
   sequence of radii. Then, $\gamma^{}_{\omega}$ admits 
   a representation of the form
$$ \gamma^{}_{\omega} \; = \; g\cdot \delta^{}_{\gG} $$
   where $\delta^{}_{\gG} = \sum_{t\in\gG}\delta_t$ is the
   Dirac comb of\/ $\gG$ and $g$ is a bounded, positive 
   definite Lipschitz function on all of\/ $\RR^n$. \qed
\end{prop}

Now, we can invoke Bochner's Theorem \cite[Thm.\ IX.9]{RS} which 
tells us that  $\hat{g}$, the Fourier transform of $g$, is a {\em finite\/} 
positive measure (our above construction was crucial to ensure this).
On the other hand, $\hat{\delta}^{}_{\gG}$ is translation
bounded, wherefore $\hat{g} * \hat{\delta}^{}_{\gG}$ is well defined
\cite[Prop.\ 1.13]{BF}, and it is a tempered measure.
This allows us to (backwards) employ the convolution theorem,
compare \cite[Prop.\ 4.10]{BF}:
\begin{equation} \label{rep2}
 \hat{\gamma}^{}_{\omega} \; = \; \hat{g} * \hat{\delta}^{}_{\gG}
   \; = \; {\rm dens}(\gG)\,\hat{g} * \delta^{}_{\gG^*} 
\end{equation}
where the second step uses the identity
$ \hat{\delta}^{}_{\gG} \; = \; {\rm dens}(\gG)\,\delta^{}_{\gG^*}$,
known as Poisson's summation formula for lattice Dirac combs, 
see \cite[Ex.\ 6.22]{BF} or \cite{Cordoba}.
\begin{prop} \label{period1}
   The diffraction measure\/ $\hat{\gamma}^{}_{\omega}$ of\/
   $(\ref{rep2})$ is\/ $\gG^*$-periodic. It can thus alternatively
   be represented as
$$ \hat{\gamma}^{}_{\omega} \; = \; \varrho * \delta^{}_{\gG^*} $$
   where $\varrho$ is a bounded positive measure that is supported on a
   fundamental domain of\/ $\gG^*$ which we may choose to be bounded.
\end{prop}
{\sc Proof}: The $\gG^*$-periodicity of $\hat{\gamma}^{}_{\omega}$
is obvious from the convolution formula in (\ref{rep2}). 
Also, it is a standard result that each lattice in $\RR^n$
has a fundamental domain inside the Voronoi cell of $\gG^*$ that
is a Borel set, $E$ say (the construction of such a set was described
in the Remark following Theorem~\ref{thm1}). 
The translates $t+E$ for $t\in\gG^*$ then
form a disjoint partition of $\RR^n$, and we have
$$\hat{\gamma}^{}_{\omega} \; = \; \sum_{t\in\gG^*}
   \hat{\gamma}^{}_{\omega}\big|^{}_{t+E} $$
where $\hat{\gamma}^{}_{\omega}\big|^{}_{t+E}$ is the restriction of
$\hat{\gamma}^{}_{\omega}$ to the set $t+E$.
But periodicity tells us that $\hat{\gamma}^{}_{\omega}\big|^{}_{t+E}
= \hat{\gamma}^{}_{\omega}\big|^{}_{E} * \delta_t$, and the
assertion follows with $\varrho:=\hat{\gamma}^{}_{\omega}\big|^{}_{E}$.
\qed

\smallskip
To complete the proof of Theorem~\ref{thm1}, it remains to go back 
to the autocorrelation, i.e.\
$$ \gamma^{}_{\omega} \; = \; \check{\varrho}\cdot
   \check{\delta}^{}_{\gG^*} \; = \; \gP \cdot \delta^{}_{\gG} $$
where we have
$$ \varrho \; = \; {\rm dens}(\gG)\, \hat{\gP} \, ,$$
again by Poisson's summation formula. The nice properties of
the function $\gP$ claimed in Theorem \ref{thm1} are now a direct
consequence of the Paley-Wiener Theorem, see \cite[Thm.\ IX.12]{RS}
or \cite[Thm.\ 7.23]{Rudin2}, because 
$\varrho$ is by construction a positive measure with compact support.
This also concludes the proof of Theorem \ref{thm1}.
   \qed

\smallskip
\noindent
{\sc Remark}: The above proof employed the Lipschitz property of the
bump functions $c$. Alternatively, using (\ref{faltung}) and the 
smallness of the support of $c$ (and hence also of that of 
$c*\tilde{c}$), one can derive that
$$ \big| g^{}_r (z) - g^{}_r (z') \big| \; \le \;
   \frac{|\gG_r|}{{\rm vol}(B_r(0))}\, W^2\,
   \sup_{x} \big|\big(c*\tilde{c}\big) (z-x) - 
                 \big(c*\tilde{c}\big) (z'-x)\big|  $$
for all $z,z'$ sufficiently close. This inequality implies equi-continuity of the 
family $\{g^{}_r\mid r>0\}$ without reference to any Lipschitz property.
This results in a slightly weaker version of Proposition \ref{rep1},
with $g$ being a bounded, positive definite continuous function on all of
$\RR^n$ which need not be Lipschitz.

\section{\Large Complementary lattice subsets}

A particularly interesting situation emerges in the comparison of
a lattice subset $S\subset\gG$ with its complement $S'=\gG\setminus S$.
Here, the Dirac combs to be compared are 
$\omega = \delta_S = \sum_{x\in S}\delta_x$
and $\omega' = \delta_{S'}$. Note that the Dirac comb of (\ref{comb1}),
when specialized to $w\equiv 1$, results in an autocorrelation
with coefficients
$$ \nu^{}_S(z) \; = \; \lim_{r\to\infty}\,
   \frac{|S \cap (z+S) \cap B_r(0)|}{{\rm vol} (B_r(0))}
   \; = \; \dens\big(S\cap(z+S)\big) $$
for all $z\in\RR^n$, provided the limits exist. This can easily be derived 
from (\ref{auto2}) and the comments following it.

Next, recall that two point sets are called {\em homometric\/} if they
share the same (natural) autocorrelation. This is an important concept
in crystallo\-graphy, both in theory and practice, because 
homometric sets cannot be distinguished by diffraction \cite{Pat,RoSe}.

\begin{theorem} \label{complement}
   Let $\gG$ be a lattice in $\RR^n$, and let
   $S\subset\gG$ be a subset with existing 
   $($natural\/$)$ autocorrelation
   coefficients $\nu^{}_{S}(z) = \dens\big(S\cap(z+S)\big)$.
   Then the following holds.

$(1)$ The autocorrelation coefficients $\nu^{}_{S'}(z)$ of the 
   complement set $S'=\gG\setminus S$ also exist. They are
   $\nu^{}_{S'}(z)=0$ for all $z\not\in\gG$
   and otherwise, for $z=t\in\gG$, satisfy the relation
$$ \nu^{}_{S'}(t) - \dens(S')  \; = \;
   \nu^{}_{S}(t) - \dens(S) \, . $$

$(2)$  If, in addition, $\dens(S) = \dens(\gG)/2$, then the sets
    $S$ and $S'=\gG\setminus S$ are homometric.

$(3)$ The diffraction spectra of the sets $S$ and $S'$ are related by
$$ \hat{\gamma}^{}_{S'} \; = \; \hat{\gamma}^{}_{S} 
   + \big(\dens(S') - \dens(S)\big)\, \dens(\gG)\, \delta^{}_{\gG^*} \, . $$
   In particular, $\hat{\gamma}^{}_{S'} = \hat{\gamma}^{}_{S}$ if\/
   $\dens(S') = \dens(S)$.

$(4)$ The diffraction measure $\hat{\gamma}^{}_{S'}$ is pure point
   if and only if $\hat{\gamma}^{}_{S}$ is pure point.
\end{theorem}
{\sc Proof}: In what follows, each term involving a density is to be
viewed as the limit along a fixed increasing and unbounded sequence of radii.
Since $\gG$ is the disjoint union of $S$ and $S'$, 
$\gG = S \, \dot{\cup} \, S'$, we get
$\dens(S') = \dens(\gG) - \dens(S)$ and the natural density of $S'$ exists
because $\dens(S) = \nu^{}_{S} (0)$. Since $S'\subset\gG$, we also have
$\nu^{}_{S'}(z)=0$ whenever $z\not\in\gG$.

So, let $z=t\in\gG$ from now on. Now observe that $\gG\cap(t+\gG) = \gG$
and thus, using $\gG = S \, \dot{\cup} \, S'$, we obtain
\begin{eqnarray*}
   \dens(\gG) & = & \dens\big(\gG\cap(t+\gG)\big) \\
   & = & \nu^{}_{S'}(t) + \nu^{}_{S}(t) + 
   \dens\big(S\cap(t+S')\big) + \dens\big(S'\cap(t+S)\big). 
\end{eqnarray*}
Since $S'=\gG\setminus S$, it is easy to verify that
\begin{eqnarray*}
   \dens\big(S'\cap(t+S)\big) & = &
   \dens\big(\gG\cap(t+S)\big) - \dens\big(S\cap(t+S)\big) \\
   & = & \dens(S) - \nu^{}_{S}(t)
\end{eqnarray*}
because $(t+S)\subset\gG$ and $\dens(t+S)=\dens(S)$. Similarly,
$$  \dens\big(S\cap(t+S')\big) \; = \;
    \dens(S) - \nu^{}_{S}(-t) $$
by first shifting and then using the previous formula.
Since $\nu^{}_{S}(t)$ is a real positive definite function, we have
$\nu^{}_{S}(-t)=\nu^{}_{S}(t)$, and obtain
$$ \dens(\gG) = 2\, \dens(S) + \nu^{}_{S'}(t) - \nu^{}_{S}(t) $$
from which the first assertion follows with 
$\dens(\gG) = \dens(S) + \dens(S')$.   

If $\dens(S) = \dens(\gG)/2$, then $\dens(S')=\dens(S)$
and we obtain $\nu^{}_{S'}(z) = \nu^{}_{S}(z)$, for all $z$, by the
first assertion. This settles assertion (2).

Since $S\subset\gG$, its autocorrelation is 
$\gamma^{}_{S} = \sum_{t\in\gG} \nu^{}_{S}(t) \delta_t$,
and analogously for $S'$, the complement set in $\gG$.
From the first assertion, we then infer
$$  \gamma^{}_{S'} \; = \; \gamma^{}_{S} +
        c\, \delta^{}_{\gG}   $$
with $c = \dens(S') - \dens(S)$. Assertion (3) now follows from
taking the Fourier transform and applying
Poisson's summation formula to the lattice Dirac comb
$\delta^{}_{\gG}$. 

{}Finally, the difference between $\hat{\gamma}^{}_{S'}$ and $\hat{\gamma}^{}_{S}$
in the third assertion is a multiple of $\delta^{}_{\gG^*}$ which
is a uniform lattice Dirac comb and hence a pure point measure, whence
the last claim is obvious. \qed

\smallskip
In \cite{BMP}, it was shown that the set of visible lattice points is
pure point diffractive. The last assertion of Theorem~\ref{complement} 
then tells us that their complement, the set of {\em invisible\/} points, 
also is pure point diffractive. Similarly, the set of $k$-th power free integers,
a subset of $\ZZ$, has pure point diffraction \cite{BMP}, so does then
its complement, the set of integers divisible by the $k$-th power of
some integer $\ge 2$. This indicates that many more pure point diffractive
point sets of independent interest exist, and a general criterion based
on the almost periodicity of the autocorrelation is derived in \cite{BMnew}.

\section{\Large Generalizations}

Our above derivation, with little modification, can also be carried
through in the case that $\RR^n$ is replaced by an arbitrary locally compact
Abelian group $G$ which is $\sigma$-compact and metrizable.  
As such, $G$ is certainly Hausdorff, but also Polish, see the Corollary 
in \cite[Ch.\ IX.6.1]{Bou}. In particular, we may
assume $G$ to be equipped with a metric $d$ which induces the topology of 
$G$ and with respect to which $G$ is complete. This class contains all groups
of the form $\RR^n\times\ZZ^m\times H$ with $n,m\ge 0$ and $H$ compact and
metrizable, which are (up to isomorphism) the so-called 
metrizable compactly generated LCA groups. We assume $G$ to be equipped with
an appropriately normalized Haar measure $\theta^{}_G$, compare 
\cite[Ch.\ 2.5]{Gu}. 
In particular, we assume $\theta^{}_{G}(G)=1$ if $G$ is a compact group.

A {\em lattice\/} $\gG$ is now a discrete (hence closed) subgroup of $G$ such that 
$G/\gG$ is compact (this is the appropriate generalization of our previous
definition in the Euclidean setting). 
The Dirac comb $\omega$ of (\ref{comb1}) is well defined, 
and it is translation bounded (or shift bounded) if 
and only if the function $w$ is bounded.
Note that $\omega$ is then a {\em finite\/} complex measure if $G$ is compact.
The (open) ball of radius $r$ around $a$ is 
$$ B_r(a) = \{ x\in G \mid d(x,a) < r \}\, , $$ 
but the autocorrelation $\gamma^{}_{\omega}$ of the Dirac 
comb $\omega$ can, in general, no longer be defined as in Eq.~(\ref{auto1}), 
because balls in general groups $G$ can have rather weird properties. So, even with
${\rm vol}(B_r(0)) := \theta^{}_{G}(B_r(0))$, the limit in (\ref{auto1}) might be
meaningless. There is no problem for compact $G$, though: one can simply write
$\gamma^{}_{\omega} = \omega * \tilde{\omega}$ because $\theta^{}_{G}(G)=1$
in this case, so that existence and uniqueness of the autocorrelation are
automatic. This is clear since all sums involved are actually finite.

We overcome the general difficulty by employing the 
concept of an {\em averaging sequence\/} which, at the same time,
also constitutes a restricted (monotone increasing) van Hove sequence, 
see \cite{Sch} for details. To explain this, recall that an LCA group $G$ is 
$\sigma$-compact if and only if a countable family
\begin{equation} \label{sequence}
   {\mathcal U} \; = \; \{ U_i \mid i\in \NN\, \}
\end{equation}
of relatively compact open sets exists with $U_1\neq \varnothing$,
$\overline{U_i}\subset U_{i+1}$ for all $i\in\NN$,
and $\bigcup_{i\in\NN} U_i = G$, see \cite[Thm.\ 8.22]{Q}.
In particular, $0<\theta^{}_{G}(U_i)<\infty$ for all $i\in\NN$.
We call such a family $\mathcal U$ an averaging sequence. It also
constitutes a {\em van Hove sequence\/}, if one extra property is satisfied
which guarantees that the surface/bulk ratio becomes sufficiently negligible
in the limit $i\to\infty$. To formalize the latter, we introduce
$$ \partial^K U \; = \; \big( (\, \overline{U} + K ) \setminus U \big)
   \cup \big( ( G\setminus U - K )\cap \overline{U}\, \big)  $$
for an open set $U$ and any compact set $K$ in $G$. Here, we use the 
convention $A\pm B := \{x\pm y \mid x\in A, y \in B \}$, with 
$\varnothing\pm B = \varnothing$.
The set $\partial^K U$ can, cum grano salis, be seen as a `thickened' 
version of the boundary of $U$. Then, the final condition is that
\begin{equation} \label{hove}
   \lim_{i\to\infty} \frac{\theta^{}_{G}(\partial^K U_i)}{\theta^{}_{G}(U_i)} 
   \; = \; 0
\end{equation}
for all compact $K\subset G$. The existence of such averaging sequences
of van Hove type in $\sigma$-compact LCA groups is shown in \cite{Sch}. 

Let us assume that a family ${\mathcal U}=\{U_i\mid i\in\NN\}$ has been 
given which also constitutes a van Hove sequence. Since $\gG$, as a lattice, 
is a special case of a regular model set, we can then conclude from 
\cite[p.~145]{Sch} that the {\em density\/} of $\gG$ exists,
\begin{equation} \label{density2}
   {\rm dens}(\gG) \; = \; \lim_{i\to\infty}\,
   \frac{|\gG\cap U_i|}{\theta^{}_{G}(U_i)} ,
\end{equation}
and that the limit is independent of the van Hove sequence chosen. This is the 
correct analogue of the statement (\ref{density}) for lattices in $\RR^n$.
If we now define $\omega^{}_i = \omega |^{}_{U_i}$ and, similarly,
$\tilde{\omega}^{}_i = (\omega |^{}_{U_i})\tilde{\hphantom{t}}$, the analogue of 
Eq.~(\ref{auto1}) would read
\begin{equation} \label{new-auto1}
    \gamma^{}_{\omega} \; = \; \lim_{i\to\infty}
    \frac{\omega^{}_i * \tilde{\omega}^{}_i}{\theta^{}_{G}(U_i)} 
\end{equation}
provided the limit exists.

In this case, the autocorrelation is always a pure point measure
of the form $\gamma^{}_{\omega} = \sum_{z\in\gG} \nu(z)\delta_z$
with the coefficients $\nu(z)$ now being given by
\begin{equation} \label{new-auto2}
    \nu(z) \; = \; \lim_{i\to\infty} \,
    \frac{1}{\theta^{}_{G}(U_i)} \, 
    \sum_{\substack{t,t'\in\gG\cap U_i \\ t-t'=z}}
    {w(t)}\, \overline{w(t')}\, , 
\end{equation}
again with the simplification that 
$\nu(z) = \sum_{t\in\gG} w(t) \overline{w(t-z)}$
for $G$ compact. Usually, the more interesting cases of Dirac combs will 
occur for groups $G$ that are locally compact, but {\em not\/} compact, 
such as $\RR^n$.

At this point, let us state the analogue of Theorem~\ref{thm1} in this more
general setting.
\begin{theorem} \label{thm2}
   Let $\gG$ be a lattice in a Polish LCA group $G$, and let
   $\omega = \sum_{t\in\gG} w(t)\, \delta_t$ be a Dirac comb on $\gG$
   with bounded complex weights.
   Let an averaging sequence ${\mathcal U} = \{U_i\mid i\in\NN\}$ be given which 
   constitutes a monotone increasing van Hove sequence, and assume
   that the corresponding autocorrelation measure $\gamma^{}_{\omega}$ of\/ 
   $(\ref{new-auto1})$ exists. Then we have:

$(1)$   The autocorrelation measure $\gamma^{}_{\omega}$ admits the representation
$$ \gamma^{}_{\omega} \; = \; g \cdot \delta^{}_{\gG} $$
   where $g$ is a bounded, positive definite Lipschitz function on\/ $G$
   which interpolates the autocorrelation coefficients $\nu(t)$, $t\in\gG$.

$(2)$   The Fourier transform $\hat{\gamma}^{}_{\omega}$ is a translation
   bounded positive measure on the dual group $\hat{G}$ which is periodic with 
   lattice of periods\/ $\gG^*$, the dual lattice of\/ $\gG$. 
   It can be represented as
$$ \hat{\gamma}^{}_{\omega} \; = \; \varrho * \delta^{}_{\gG^*} $$
   with a finite positive measure $\varrho$ that is supported on a 
   totally bounded and measurable fundamental domain of\/ $\gG^*$.  
\end{theorem}

{\sc Proof}: The argument is very similar to the proof in Section~\ref{proofs},
wherefore we only describe the changes needed.

Lemma \ref{lip} needs a replacement because the concept of a smooth bump
function makes no sense in general. Instead, we can employ a different kind
of Lipschitz function as follows. Choose $\varepsilon > 0$ so that 
$B_{\varepsilon}(0)\cap B_{\varepsilon}(t) =\varnothing$ for all 
$t\in\gG\setminus\{0\}$. Such an $\varepsilon$ clearly exist because
$G$ is Hausdorff and our metric induces the topology of $G$. If we set 
$A=G\setminus B_{\varepsilon}(0)$, which is a closed set, then $0\not\in A$
and thus $d(0,A)\ge\varepsilon>0$, where $d(x,A) = \inf_{y\in A} d(x,y)$ is 
the distance of $x$ from $A$. We can now define 
\begin{equation} \label{new-lip}
     c(x) \; = \; \frac{c^{}_0 \, d(x,A)}{d(0,A)}\, . 
\end{equation}
One has $c(x)\ge 0$, $c(0)=c^{}_0$ and $c(x)=0$ for all $x\in A$, so $c$ 
is nontrivial
(if $c^{}_0\neq 0$) and supported on $B_{\varepsilon}(0)$.
\begin{lemma} \label{lemma1a}
   The function $c$ of $(\ref{new-lip})$ is Lipschitz, with Lipschitz
   constant $L_c = c^{}_0/d(0,A)$. The functions $\tilde{c}$ and
   $c*\tilde{c}$ als also Lipschitz, with $L_{\tilde{c}} = L_c$
   and $L_{c*\tilde{c}}\le \|c\|^{}_{1} L_c$.
\end{lemma}
{\sc Proof}: The Lipschitz property of $d(x,A)$, with Lipschitz constant $1$,
is stated in \cite[Prop.\ IX.2.3]{Bou}, see \cite{BM-book} for an explicit
proof. The remainder of the proof is that of Lemma~\ref{g-lip}, with
$\RR^n$ and Lebesgue measure replaced by $g$ and Haar measure, respectively.
\qed

\smallskip
Since both $c$ and $c*\tilde{c}$ are integrable w.r.t.\ the Haar measure 
$\theta^{}_{G}$, the construction 
of a family of approximating functions $\{ g^{}_i \mid i\in\NN \}$,
relative to the sequence $\mathcal U$, is now possible,
in complete analogy to above. 

With this modification, Lemma \ref{g-lip} still holds. For its proof, we only 
need that the density of lattice points of $\gG$ in $G$ exists and that, if
$G$ is not compact, $|\gG\cap U_i| / \theta^{}_{G}(U_i)$ converges to it 
as $i\to\infty$. This follows from (\ref{density2}).
We can then proceed as before: since we assume $G$ to be 
$\sigma$-compact, we have Ascoli's Theorem at our disposal. Consequently, 
we obtain the following modification of Proposition \ref{rep1}.
\begin{prop} \label{lca-rep1}
   Let $\gG$ be a lattice in a Polish LCA group $G$ and
   let $\omega$ be the Dirac comb of $(\ref{comb1})$, with 
   bounded complex weights. Let $\gamma^{}_{\omega}$ be an
   autocorrelation measure which is assumed to exist as the limit 
   of a sequence $\{\gamma^{(i)}_{\omega} \mid i\in\NN\}$  of
   approximants with respect to a given averaging sequence\/ 
   ${\cal U} = \{U_i\mid i\in\NN\}$. Then, $\gamma^{}_{\omega}$ admits 
   a representation of the form
$$ \gamma^{}_{\omega} \; = \; g\cdot \delta^{}_{\gG} $$
   where  $g$ is a bounded, positive 
   definite Lipschitz function on all of\/ $G$. \qed
\end{prop}

Next, we need some results from harmonic analysis in the setting of LCA groups,
see \cite[Ch.\ I.4]{BF} or \cite[Ch.\ 2.8]{Gu} for a suitable summary. We 
denote the dual group
by $\hat{G}$ and assume it is equipped with a matching Haar measure,
$\theta^{}_{\hat{G}}$ (a suitable normalization of it is suggested by the Fourier
inversion formula, see \cite[Ch.\ 2.8.7]{Gu}).
We can now, once more, invoke Bochner's Theorem \cite[Thm.\ 4.5]{BF}:
since $\gamma^{}_{\omega}$ is a {\em positive definite\/} measure by construction,
its Fourier transform is a uniquely determined translation bounded 
{\em positive\/} measure on $\hat{G}$, denoted by $\hat{\gamma}^{}_{\omega}$. 

To make complete sense out of Eq.~(\ref{rep2}),
we have to say what the dual lattice is and how Poisson's summation formula
works. Each $k\in\hat{G}$ defines a (continuous) character on the group
$G$, $\langle k,x\rangle$, which replaces $\exp(2\pi i kx)$ from above. Then,
$$  \gG^* \; = \; \{ k \in \hat{G} \mid \langle k,x\rangle = 1
    \mbox{ for all } x \in \gG \}  $$
and $\gG^*$ (which is called $\gG^{\perp}$ in \cite{Gu}) is the 
{\em annihilator\/} 
of the closed subgroup $\gG\subset G$ in the dual group $\hat{G}$, compare 
\cite[Ch.\ 2.9.1]{Gu}. We have $\gG^* = (G/\gG)\!\hat{\hphantom{A}}$ and
$\hat{\gG} = \hat{G}/\gG^*$, see \cite[Thm.\ 2.9.1]{Gu}, so that
with $\gG$ also $\gG^*$ is a lattice because the dual of a compact
group is discrete and vice versa \cite[Thm.\ 2.8.3]{Gu}.
Moreover, we can interpret $G/\gG$ and $\hat{G}/\gG^*$ as 
(measurable) fundamental
domains of the lattices $\gG$ and $\gG^*$. We now get, in complete analogy 
to before, the general Poisson summation formula
\begin{equation}  \label{gen-poisson}
    \hat{\delta}^{}_{\gG} \; = \; a \, \delta^{}_{\gG^*}
\end{equation}
where $a$ is a constant which depends on the density of $\gG$ and
on the relative normalization
of $\theta^{}_{G}$ and $\theta_{\hat{G}}$, see \cite[Ch.\ 9.9]{Gu} for details.
So, the following modification of Proposition \ref{period1} holds in our more 
general setting.

\begin{prop} \label{lca-period1}
   Under the assumptions of Proposition~$\ref{lca-rep1}$,
   the corresponding diffraction measure\/ $\hat{\gamma}^{}_{\omega}$
   exists and is\/ $\gG^*$-periodic. It can  
   be represented as
$$ \hat{\gamma}^{}_{\omega} \; = \; \varrho * \delta^{}_{\gG^*} $$
   where $\varrho$ is a bounded positive measure, supported on a
   fundamental domain of\/ $\gG^*$ which we may choose to be bounded.
   \qed
\end{prop}

Finally, we have to combine Propositions \ref{lca-rep1} and \ref{lca-period1} 
which completes the proof of Theorem~\ref{thm2}. \qed

\smallskip
\noindent
{\sc Remark}:
Similar to the situation in Theorem~\ref{thm1}, the assumption on the
existence of the autocorrelation as a limit is not essential. If
$\gamma^{(i)}_{\omega}:=\omega^{}_i * \tilde{\omega}^{}_i / \theta^{}_{G}(U_i)$,
the uniform translation boundedness of the $\gamma^{(i)}_{\omega}$
implies the existence of at least one limit point, $\gamma$ say. We can 
then choose a subfamily $\{U_{i_j}\mid j\in \NN\}$ such that 
$\gamma^{(i_j)}_{\omega} \to \gamma$ vaguely, as $j\to\infty$, and
we thus obtain the corresponding results for each 
vague limit point $\gamma$ of $\{\gamma^{(i)}_{\omega} \mid i\in\NN \}$ 
separately, by applying Theorem~\ref{thm2} to $\omega$ together with
this averaging subsequence.
 
\smallskip
{}Finally, if one needs even more generality, one can employ the Remark
at the end of Section~\ref{proofs}, which remains valid in this more
abstract setting if elements of $\cal U$ are used instead of balls, 
and start with a suitable continuous function
$c$ of sufficiently small support. If $G$ is an LCA group (which includes $G$
being Hausdorff in our terminology) and $\sigma$-compact, we obtain
\begin{theorem} \label{thm4}
   Let $G$ be a $\sigma$-compact LCA group
   $($not necessarily Polish\/$)$,
   and let the other assumptions
   be as in Theorem~$\ref{thm2}$. Then, the first statement of 
   Theorem~$\ref{thm2}$ is still true after the modification that $g$
   is merely continuous $($not necessarily Lipschitz\/$)$, while the second
   statement remains unaltered.  \qed
\end{theorem}

\bigskip

\subsection*{Acknowledgments}

It is my pleasure to thank Robert V.\ Moody, Peter A.\ B.\ Pleasants,
and Martin Schlottmann for a number of clarifying discussions, 
and the Department of Mathematical Sciences of the University
of Alberta (Edmonton, Canada) for hospitality during a stay in summer 2000,
where part of this work was done. I also thank an anonymous referee for
a number of suggestions to improve the results.
It was also supported by the German Research Council (DFG).

\end{document}